\documentclass [12pt]{article}

\begin{document} 
\newtheorem{example}{Example}[section]
\newtheorem{theorem}{Theorem}[section]
\newtheorem{lemma}{Lemma}[section]
\newtheorem{corollary}{Corollary}[section]
\newtheorem{proposition}{Proposition}[section]
\newtheorem{remark}{Remark}[section]
\setlength{\textwidth}{13cm}
\def \halmos{\hfill\mbox{ qed}\\} 
\newcommand{\eqnsection}
{\renewcommand{\theequation}{\thesection.\arabic{equation}}
\makeatletter \csname  @addtoreset\endcsname{equation}{section}
\makeatother}
 
%

\def \njc{{\bf !!! note Jay's change  !!!}}
\def  \enc{{\bf end of current changes }}
\def \nnc{{\bf !!! note NEW	 change  !!!}}
\def\square{{\vcenter{\vbox{\hrule height.3pt
                    \hbox{\vrule width.3pt height5pt \kern5pt
                       \vrule width.3pt}
                    \hrule height.3pt}}}}
\def \grad{\bigtriangledown}
\def \nc{{\bf !!! note change !!! }}
        \def \Proof{\noindent{\bf Proof}$\quad$}
\newcommand{\re}[1]{(\ref{#1})}
\def \ov{\overline}
\def \un{\underline}
\def \be{\begin{equation}} 
\def \ee{\end{equation}}
\def \bt{\begin{theorem}}
\def \et{\end{theorem}}
\def \bc{\begin{corollary}}
\def \ec{\end{corollary}}
\def \br{\begin{remark} }
\def \er{ \end{remark}}
\def \bl{\begin{lemma}}
\def \el{\end{lemma}}
\def \bex{\begin{example}}
\def \eex{\end{example}}
\def \bea{\begin{eqnarray}}
\def \eea{\end{eqnarray}}
\def \bas{\begin{eqnarray*}}
\def \eas{\end{eqnarray*}}
\def \al{\alpha}
\def \bb{\beta}
\def \ga{\gamma}
\def \Ga{\Gamma}
\def \de{\delta}
\def \De{\Delta}
\def \ep{\epsilon}
\def \vep{\varepsilon}
\def \la{\lambda}
\def \La{\Lambda}
\def \ka{\kappa}
\def \om{\omega}
\def \Om{\Omega}
\def \va{\varrho}
\def \ffi{\Phi}
\def \vf{\varphi}
\def \si{\sigma}
\def \Si{\Sigma}
\def \vsi{\varsigma}
\def \th{\theta}
\def \Th{\Theta}
\def \ups{\Upsilon}
\def \ze{\zeta}
\def \tr{\nabla}
\def \ff{\infty}
\def \wh{\widehat}
\def \wt{\widetilde}
\def \dar{\downarrow}
\def \rar{\rightarrow}
\def \uar{\uparrow}
\def \sbs{\subseteq}
\def \mpt{\mapsto}
\def \R{{\bf R}}
\def \G{{\bf G}}
\def \H{{\bf H}}
\def \Z{{\bf Z}}
\def \S{{\bf
S}}
\def \sfB{{\sf B}}
\def \sfS{{\sf S}}
\def \T{{\bf T}}
\def
\C{{\bf C}}
\def \AA{{\mathcal A}}
\def \BB{{\mathcal B}}
\def
\CC{{\mathcal C}}
\def \DD{{\mathcal D}}
\def \EE{{\mathcal E}}
\def
\FF{{\mathcal F}}
\def \GG{{\mathcal G}}
\def \HH{{\mathcal H}}
\def
\II{{\mathcal I}}
\def \JJ{{\mathcal J}}
\def \KK{{\mathcal K}}
\def
\LL{{\mathcal L}}
\def \MM{{\mathcal M}}
\def \NN{{\mathcal N}}
\def
\OO{{\mathcal O}}
\def \PP{{\mathcal P}}
\def \QQ{{\mathcal Q}}
\def
\RR{{\mathcal R}}
\def \SS{{\mathcal S}}
\def \TT{{\mathcal T}}
\def
\UU{{\mathcal U}}
\def \VV{{\mathcal V}}
\def \ZZ{{\mathcal Z}}
\def \Pxh{P^{x/h}}
\def \Exh{E^{x/h}}
\def \Px{P^{x}}
\def \Ex{E^{x}}
\def
\Prh{P^{\rho/h}}
\def \Erh{E^{\rho/h}}
\def \p{p_{t}(x,y)}
\def
\({\left(}
\def \){\right)}
\def \lk{\left[}
\def \rk{\right]}
\def
\lc{\left\{}
\def \rc{\right\}}
\def \bsq{\ $\Box$}
\def
\nn{\nonumber}
\def \Bo{\bigotimes}
\def \bo{\times}
\def
\ot{\times}
\def
\bs{\begin{slide} }
\def \es{\end{slide} }
\def \bpr{\begin{proof}
}
\def \epr{\end{proof} }
\def \cd{\,\cdot\,}
\def
\st{\stackrel{def}{=}}
\def \as{almost surely }
\def \fix {{\bf
!!!!! }}
\def \stl{\stackrel{law}{=}}
\def \std{\stackrel{dist}{=}}
\def \stdto{\stackrel{dist}{\longrightarrow}}
\def  \enc{{ \bf end of current changes}}
\def\square{{\vcenter{\vbox{\hrule height.3pt
                    \hbox{\vrule width.3pt height5pt \kern5pt
                       \vrule width.3pt}
                    \hrule height.3pt}}}}
\def\qed{{\hfill $\square$ }}

           \def \tb{|\!|\!|}

        \eqnsection
\bibliographystyle{amsplain}

\title{  Asymptotic expansions for functions of the increments of   certain  Gaussian
processes}
        \author{ Michael B. Marcus\,\, Jay Rosen \thanks{Research of 
both authors
supported by  grants from the National Science Foundation and PSC-CUNY.}}


\maketitle
\eqnsection

\bibliographystyle{amsplain}

          \def \tb{|\!|\!|}

\begin{abstract}    Let  $G=\{G(x),x\ge 0\}$ be a mean zero
Gaussian process with stationary increments  and set $\sigma ^2(|x-y|)= 
E(G(x)-G(y))^2$.     
Let  $f$ be a    function with 
$Ef^{2}(\eta)<\ff$, where $\eta=N(0,1)$. 
When  $\si^2$
is   regularly varying at zero and   
\[
 \lim_{h\to 0}{h^2\over \sigma^2(h)}= 0\qquad \mbox{and}\qquad \lim_{h\to
0}{
\sigma^2(h)\over h}= 0
  \quad \mbox{but} \quad  \({d^{2}\over
ds^2}\si^2(s)\)^{j_0}    
\]
 is locally integrable for some integer $j_0\ge 1$, and satisfies some additional regularity
conditions, 
\bea
&& \int_a^bf\(\frac{G(x+h)-G(x)}{\sigma (h)}
\)\,dx \label{abst}\nn
 \\
&&\qquad  = \sum_{j=0}^{j_0} (h/\sigma(h))^{j}\,\,{E(H_{j}(\eta)
f(\eta))\over\sqrt {j!}}\,\,:(G')^{j}:(I_{[a,b]})\, +o\({h\over\sigma (h)}\)^{j_0}\nn
\eea
in  $L^2$. Here $H_j$ is the $j$-th Hermite
polynomial. 
Also $:(G')^{j }:(I_{[a,b]})$ is a $j $-th order Wick power Gaussian
chaos constructed from the Gaussian field
  $ G'(g) $, with covariance
\[
E(G'(g)G'(\wt g)) = \int \!\!\int  \rho (x-y)g(x)\wt g(y)
 \,dx\,dy\label{3.7bqs} ,
\]
where   $   \rho(s)=\frac{1}{2}{d^{2}\over ds^2}\sigma^2(s)$.

\end{abstract}

Ê ÊÊ Ê\section{Introduction}\label{sec-intro}

LetÊ $G=\{G(x),x\in R_+ \}$,Ê Ê $G(0)=0$, be a mean zeroÊ GaussianÊ process
with stationary increments, and set 
\be 
Ê Ê ÊÊ E(G(x)-G(y))^2=\si ^2( x-y ) \label{m1.2}=\si
^2(|x-y|).\nn
\ee 
ÊÊ The function $\si^2$ is referred to as the
increment's variance of $G$.Ê Ê Ê Clearly $\si^2(0)=0$. 

  In this paper we are
primarily concerned with Gaussian processes that are smoother than Brownian
motion but not so smooth that they have mean square
derivatives.  

  Let $d\mu( x)=( 2\pi)^{ -1/2}\exp ( -x^{ 2}/2)\,dx$
denote standard Gaussian measure on
$R^{ 1}$.  Let
$f\in L^{ 2}(R^{ 1},\,d\mu )$,Ê i.e.,   $Ef^{2}(\eta)<\ff$, where $\eta $ is a
normal random variable with mean zero and variance one, (i.e.  
$\eta= N(0,1)$).  To avoid trivialities we assume that
$\si^2(h)\not\equiv 0$ and $f(x)\not\equiv 0$.  In all that follows $0\le
a<b<\ff.$ 

\medskip

  We obtain   an $L^{2}$ asymptotic expansion  for 
\begin{equation}
 \int_a^bf\(\frac{G(x+h)-G(x)}{\si (h)}
\)\,dx\label{add.1},
\end{equation}
 as $h\rar 0$,  that holds for a large class of Gaussian processes   and for all $f\in L^{ 2}(R^{ 1},\,d\mu )$. The asymptotic expansion   involves a generalized derivative $G'$ of the
Gaussian process $G$.  

We impose the following conditions on the Gaussian processes considered here: 
 \be
\mbox{$\si^2(h)$ is regularly varying at zero of index $1\leq\bb\leq 2$;}\label{rv}
\ee
\be\lim_{h\to 0}{h^2\over \si^2(h)}= 0\qquad \mbox{and}\qquad \lim_{h\to
0}{
\si^2(h)\over h}= 0;\label{mm5.17o}
\ee
\be 
{| \si^2(s+h)+ \si^2(s-h)-2\si^2(s)|\over h^2}\le C{\si^2(s)\over
s^2}\quad\mbox{for } h\le { s\over 8 };\label{mm5.9}
\eeÊ
 \be
\mbox{$\si^2(s)$ has a second derivative for each $s\ne	0$.}\label{rv2}
\ee
  Set  
\be \rho (s):=\frac12{d^{2}\over ds^2}\si^2(s), \hspace{.2 in}s\ne	0.\label{rho}
\eeÊ
  It follows from  (\ref{mm5.17o}) that  
\be
 {d
\si^2(0)\over ds}= 0\qquad \mbox{and} \qquad\rho(0):=\lim_{h\to 0 } { 
\si^2(h)\over h^{2}}=\ff.\label{mm5.17oj} 
\ee

The next theorem is the main result in this paper.

Ê \bt \label{theo-3.2n} Let $f\in L^{ 2}(R^{ 1},\,d\mu)$ and letÊ
$G=\{G(x),x\in R_+\}$, $G(0)=0$, be a mean zero Gaussian processÊ with
stationary incrementsÊ   satisfying (\ref{rv})--(\ref{rv2}),Ê andÊ assume thatÊ
there exists a
$  \ze>0  $  such that for all  
$0<M<\ff$ we can find $C_{M}<\ff$ with  
\begin{equation}
|\rho( x)|\leq {C_{M} \over |x|^{\ze}}  := C_{M}\varphi(x),\quad\hspace{.2 in} |x|\leq
M
\label{nor14.5a}
\end{equation}
and  
\begin{equation}
 |\rho( x+h) -\rho( x)|\leq C_{M}{|h|Ê \over |x| }\,|\rho( x)|,
\hspace{.2 in}4|h|\leq |x|\leq M. \label{nor14.5}
\end{equation}
 Then  for all integers $j_{0}$, such that Ê$j_{0}\ze<1$,Ê and for all  for $b\ge a$, 
\begin{eqnarray}
&& \int_a^bf\(\frac{G(x+h)-G(x)}{\si (h)}
\)\,dx
\label{mm5.22n}\\
&&\qquadÊ = \sum_{j=0}^{j_0} (h/\si(h))^{j}\,\,{E(H_{j}(\eta)
f(\eta))\over\sqrt {j!}}\,\,:(G')^{j}:(I_{[a,b]})\, +o\({h\over\si (h)}\)^{j_0}\nn
\end{eqnarray}
in $L^2$.Ê
 \et

   There are many terms in (\ref{mm5.22n}) that require   definition. Ê Ê The functions  $\{H_{ k}( x)\}_{k=0}^\ff$ are the
\label{pageII}Hermite polynomials. The process
 $G' =\{G'(f ),f\in\BB_{0}(R_{+})\}$ is a mean zero Gaussian field  with 
\begin{equation}
 E\(G'(f)G'(\wt f)\)=\int\int\rho (t-s)\,f( s)\,\wt f( t)   \,ds\,dt\qquad\forall  f,\wt f\in\BB_{0}(R_{+}) \label{new.5d}
\end{equation}
where $\mathcal{B}_{0}(R_+)$ is  the set of bounded Lebesgue measurable functions on $R_+$ with compact support. We construct $G'$ in Section \ref{sec-gd}.  (We use the notation $G'$
becauseÊ it is a generalized derivative of the Gaussian process
$G$. This is also explained in Section \ref{sec-gd}.)

The random variable
$:(G')^{k_0}:(I_{[a,b]})$ is theÊÊ `value' ofÊÊ the $k_0$-th order Wick power
Gaussian chaos process
$ \{:(G')^{k_0}:(g),g\in\BB_0(R_+)\}$, at
$g=I_{[a,b]}$.Ê This process  is constructed from
$G'$ in Section
\ref{sec-nonnor4} and has second moment 
\be
E\(:(G')^{k_0}:(g)\)^2=k_0\,!\int \!\!\intÊ \rho^{k_0}(x-y)g(x)g(y)
 \,dx\,dy\label{3.7bq} .
\ee
  It is well known that $:(G')^{k_0}:(I_{[a,b]})$ can also be expressed as a multiple   Wiener-It\^{o} integral. We discuss this in Section \ref{sec-nonnor4}.
 
The $k$-th order Wick power of a mean zero Gaussian random variable $X$ is
\be
:X^k:\,=\sum_{ j=0}^{[ k/2]}( -1)^{ j}{Ê k\choose 2 j}E( X^{ 2 j})\,\,X^{ÊÊ k-2j
}.\label{defin}
\ee
When $X=N(0,1)$,Ê $ :X^{ k}:\, =\sqrt{k!\,}H_{k}( X)$. Therefore  

\be
:X^k:\,=\sqrt{k!\,}\si_X^kH_{k}\( {X\over \si _X}\).\label{0.9}
\ee
  $\si^2_X$ denotes the variance of $X$. 
We show in Theorem \ref{theo-9.1} that when ÊÊ 
 \be 
\({d^{2}\over
ds^2}\si^2(s)\)^{k_0}\qquad\mbox{ is locally integrable}\label{1.15}
\ee
 and satisfies an additional very mild
regularity condition Êthen
\be
\lim_{h\downarrow 0}\int
:\({G( x+h)- G( x) \over h}\)^{ k_0}:g(x) \,dx=\,:(G')^{k_0}:(g)\label{sss}
\ee
in $L^2$.   

 When   $  \rho (0) <\ff$, $G$ has a mean square derivative and
one would expect (\ref{sss}) to hold with $G'$ being the mean square derivative. Theorem \ref{theo-3.2n} shows that this  holds for all   $f\in L^{ 2}(R^{ 1},\,d\mu)$ and for a much more
general    class of Gaussian processes. 

\medskip	 The class of Gaussian processes satisfying  the hypotheses of Theorem \ref{theo-3.2n} is   very rich. This is illustrated in the next proposition. 
\begin{proposition}\label{prop-conv} Let $h$ be any 
function that is regularly varying at infinity 
with negative
index 
or is slowly varying at infinity and decreasing.   Then, for any 
$1<\bb<2$, there
exists a   Gaussian process with stationary increments for which the increments variance 
$\si^2(x)$  satisfies the hypotheses of Theorem \ref{theo-3.2n} and is such that
\be\label{sam.019}
\si^2(x)\sim|x|^{\bb }h(\log 
1/|x|)\qquad \mbox{as}\qquad x\to 0.
\ee 
\end{proposition}

   Other examples are given in  Section \ref{sec-5}.

 \medskip  For any function
$f\in L^{ 2}(R^{ 1},\,d\mu )$,
\begin{equation} f( x)=\sum_{ k=0}^{ \ff}a_{k}H_{k}(
x)\hspace{ .2in}\mbox{ inÊÊ }\,\,L^{ 2}(R^{ 1},\,d\mu ),
\label{a18.1}
\end{equation} where
\be a_{k}=\int f( x)H_{k}( x)\,d\mu( x)=EH_k(\eta)f(\eta)\label{am.m1}
\ee
Ê Ê ÊÊ and
\begin{equation}
\sum_{ k=0}^{ \ff}a_{k}^{ 2}=Ê \int |f( x)|^{ 2}\,d\mu( x)Ê <\ff.
\label{a18.2}
\end{equation} 
For a given $f\in L^{ 2}(R^{ 1},\,d\mu )$ let  
\be
Ê k_{ 0}:=Ê k_0(f)=\inf_{k\ge 1}\{ k|a_{k} \ne 0 \}.\label{b.mm.3.1}
\ee
  The integer $ k_{ 0}$ is known as the Hermite rank of $f$.

\medskip	 We have the following  corollary of  Theorem \ref{theo-3.2n}:
 \begin{corollary} \label{theo-3.2} For a givenÊ $f\in L^{ 2}(R^{ 1},\,d\mu)$Ê let $k_0$ be as
defined in (\ref{b.mm.3.1})    and let $G=\{G(x),x\in R_+\}$, $G(0)=0$, be a mean zero Gaussian processÊ with
stationary incrementsÊ   satisfying (\ref{rv})--(\ref{rv2}). 
  Assume that (\ref{1.15}) holds with $k=k_{0}$    and that 
\be
h=o\( h^2/\si^2(h) \)^{k_0}.\label{regularity}
\ee
 Then,Ê Ê for $b\ge a$
\be
Ê \lim_{h\downarrow 0}{ \int_a^bf\(\frac{G(x+h)-G(x)}{\si (h)}
\)\,dx - (b-a)Ef(\eta)\over (h/\si(h))^{k_0}}= {E(H_{k_0}(\eta)
f(\eta))\over\sqrt {k_0!}}:(G')^{k_0}:(I_{[a,b]})\label{mm5.22}
\ee
in $L^2$.Ê
 \end{corollary}
 
  (Note that (\ref{regularity}) is implied by (\ref{nor14.5a}) if $k_{0}\ze<1$).

\begin{example} {\rm  It follows immediately from (\ref{0.9}) and  (\ref{mm5.22}) that  
\begin{equation}
   \lim_{h\downarrow 0}\int_a^b:\(\frac{G(x+h)-G(x)}{h} 
\)^{k}:\,dx =\frac{1}{\sqrt{k!} }:(G')^{k }:(I_{[a,b]})\label{1.25}
   \end{equation} 
in $L^{2}$.   Remarkably, we show in    \cite{as}, that under the hypotheses of Theorem \ref{theo-3.2n},  the limit in (\ref{1.25}) is also almost sure.   
}\end{example}
  It isÊÊ clear from (\ref{sss}) that when
$k_0>1$, the limit in (\ref{mm5.22}) is not a normal random variable. We do get a normal limit  when  $k_{0}=1$, as we state in the next corollary of Theorem \ref{theo-3.2n}.

 \begin{corollary} \label{theo-3.1}   Let  Ê $f\in L^{ 2}(R^{ 1},\,d\mu)$ Êbe such that $Ê E(\eta f(\eta))\ne 0$. LetÊ   $G=\{G(x),x\in R_+\}$, $G(0)=0$, be a mean zero Gaussian processÊ with
stationary incrementsÊ   satisfying (\ref{rv})--(\ref{rv2}). Assume    $\rho(s)$ is
locally integrable.
 Then
\be
Ê \lim_{h\downarrow 0}{ \int_a^bf\(\frac{G(x+h)-G(x)}{\si (h)}
\)\,dx - (b-a)Ef(\eta)\over h/\si(h)}= (E(\eta f(\eta)))
(G(b)-G(a) )\label{mm5.10}
\ee
in $L^2$.
 \end{corollary}

  It is interesting to  compare Corollary \ref{theo-3.1}  with    the  normal central limit theorem obtained in \cite[Theorem 1.1]{clt} that holds for all Gaussian processes with concave increment's variance   and for some  Gaussian processes with convex increment's variance      but where (\ref{1.15}) does not hold for $k_{0}=2$. 

\bt\label{theo-GCLT}{\cite[Theorem 1.1]{clt} }
Assume that  $\si^2(h)$ is
concave or that $\si^{2}(h)=h^{r}$, $1<r\le 3/2$. Then for all symmetric
functions $f\in L^{ 2}(R^{ 1},\,d\mu)$Ê
\be
\lim_{h\downarrow 0} { \int_a^bf\(\frac{G(x+h)-G(x)}{\si (h)}
\)\,dx-(b-a)E f(\eta) \over
\Phi(h)}
\stackrel{law}{=} N( 0,1),Ê \label{j.1}
\ee
where $\Phi^2(h)$ isÊ the variance of the numerator.
\et

Ê    Theorem
\ref{theo-GCLT} appears similar to Corollary \ref{theo-3.1}.ÊÊ In fact     under the conditions of  Corollary \ref{theo-3.1}
 \be
 \Phi(h)\sim h\si(b-a)E(\eta f(\eta))/\si(h)
 \ee 
 as $h\to 0$.  
 However, there are important differences
between these results. Theorem \ref{theo-GCLT} applies to symmetric
functions $f$ whereas in Ê Corollary \ref{theo-3.1} we require that $
E(\eta f(\eta))\ne 0$, which excludes symmetric
functions $f$.  Indeed we see from Corollary    \ref{theo-3.2}  that if $f$ is symmetric and $E(\eta^{2} f(\eta))\ne 0$   the dominant term on the right in (\ref{mm5.22}) is
\be
{E(H_{2}(\eta)
f(\eta))\over\sqrt {2 }}:(G')^{2}:(I_{[a,b]})\label{1.26},
\ee
as long as (\ref{1.15}) holds with $k_{0}=2$. The hypotheses of  Theorem \ref{theo-GCLT} excludes processes for which  (\ref{1.15}) holds with $k_{0}=2$. It is clear that the integrability of powers of $\rho $  at the origin play a critical role in whether or not we get normal central limit theorems.

Also note that 
Ê inÊ Ê Corollary \ref{theo-3.1}, we have
convergence in
$L^2$. (See  Remark
\ref{almostsure} for further discussion along this line.) 

 \medskip  We use $f\sim g$ at zero to indicate that
$\lim_{h\downarrow 0} f(h)/g(h)=1$ and  $f\approx g$  at zero to indicate that there exist $0<C_{1}\le C_{2}<\ff$ such that   $\liminf_{h\downarrow 0}    f(h)/\newline g(h)  \ge C_{1}$ and 
$\limsup_{h\downarrow 0}  f(h)/g(h)\le C_{2}$.

\subsection{Motivation}  The motivation for this paper comes from our work \cite{lp} on the local times $ \{L^x_t,(t,x)\in R_+\times  R\}$ of the 
 real valued  symmetric L\'evy process 
 $X=\{X(t), t\in R_+\}$ with
 characteristic
function
$ E e^{i\lambda X(t)} = e^{-t\psi(\lambda)} $. We show that if 
\[
     \si^2_0(x) \label{mm0.8}=\frac4\pi\int\limits ^\infty_0\,
       \frac{   \sin^2  \frac{\lambda x}{2}}{
     { \psi(\lambda)}}\,\,d\la 
\]    is concave, and satisfies some additional very weak
regularity conditions,
  then for any  $  p\ge 1$, and all $t\in R_+$
\[
\lim_{ h\downarrow 0}  \int_{a}^{ b} \bigg|{  L^{ x+h}_{ t} -L^{ x }_{
t}\over\si_0(h)}\bigg|^p\,dx =2^{p/2}E|\eta|^p
\int_a^b |L^{ x }_{ t}|^{ p/2}\,dx
\] for all
$a,b
$ in the extended real line  almost surely, and also  in $L^m$, $m\ge 1$.

  This result is obtained via the Eisenbaum Isomorphism Theorem and
depends on a related result for Gaussian processes $\{G(x),x\in R^1\}$ with stationary
increments. If the increments variance $\si_0^2(x)$ is concave, and satisfies some additional very weak
regularity conditions    we show in \cite{lp} that ,
\be
\lim_{h\to 0} \int_a^b\bigg|\frac{G(x+h)-G(x)}{\si_0(h)}\bigg|^p\,dx
=E|\eta |^p (b-a)\label{lp.0}
\ee for  all
$a,b\in R^1$, almost surely.  Viewing this as a strong law we then obtained the corresponding central limit theorem,  \cite[Theorem 1.1]{clt}, which we repeat as Theorem \ref{theo-GCLT} in this paper.   Initially, the motivation for this paper was to see what happens for  
   Gaussian processes that are smoother than those that satisfy the hypotheses    of Theorem
\ref{theo-GCLT}, but are not so smooth that they are mean square differentiable.  However, now that we have the results of  \cite{lp,clt} and this paper, we have an overview that enables us to present this work as method for finding limits of a natural sequence of stationary Gaussian processes.

 \medskip	 
 Let $G$ be the Gaussian process with stationary increments introduced at the very beginning this  section. Since
 \bea
 &&E\({(G(x+h)-G(x))(G(y +h)-G(y))\over \si^{2}(h )}\)\label{6.1}\\
 &&\qquad= { \si^2(x-y+h)+ \si^2(x-y-h)-2\si^2(x-y)\over \si^2(h)}\nn
   \eea
   we see that 
   \begin{equation}
  \GG_{h}(x)\st
   {G(x+h)-G(x) \over \si (h )} \qquad\qquad x\in R^{1}\label{6.2}
   \end{equation}
   is a stationary Gaussian process with $E(  \GG^{2}_{h}(0))=1$. A natural question is to ask whether 
   \begin{equation}
    \GG_{0}(x)\st \lim_{h\to 0} \GG_{h}(x)\label{6.3}
   \end{equation}
   exists. (The natural limit would be in $L^{2}$.) A necessary condition for such a limit is that the limit of the covariance  $E( \GG_{h}(x) \GG_{h}(y))$ should exist. This is given in (\ref{6.1}) which we write as 
  \be 
  E( \GG_{h}(x) \GG_{h}(y))\label{6.4} 
 = { \si^2(x-y+h)+ \si^2(x-y-h)-2\si^2(x-y)\over h^{2}}{h^{2}\over  \si^2(h)}.
   \ee   
  When $x-y\ne 0$ and $\si^{2}(s)$ has a second derivative for $s\ne 0$ 
     \be 
 \lim_{h\to 0} E( \GG_{h}(x) \GG_{h}(y))\label{6.4w} 
 =  ( \si^2)''(x-y )\lim _{h\to 0}{h^{2}\over  \si^2(h)}.
   \ee  
 (Note that even when $\si^{2}(s)$ is not differentiable at zero most Gaussian processes that one can think of have the property that  $\si^{2}(s)$ has a second derivative for $s\ne 0$. For example $\si^{2}(s)=|s|^{r}$, $0<r\le 2$ or $\si^{2}(s)=(\log 1/|s|)^{-r}\wedge 1 $, $ r\in R_{+}$.)   Thus for (\ref{6.3}) to hold $\si^{2}(s)$ must also have  a second derivative at $s= 0$.  
 
 \medskip	
 When $\si^{2}(s)$ does not  have  a second derivative at $s= 0$ we   consider a weak limit for  (\ref{6.2}),
 \begin{equation}
 \lim_{h\to 0}   \int_a^bf\(\frac{G(x+h)-G(x)}{\si (h)}
\)\,dx,\label{6.6}
   \end{equation} 
for $f\in L^{ 2}(R^{ 1},\,d\mu ) $. It is natural to approach this by first taking $f= H_{ k}( x) $ the $k$--th
 Hermite polynomial. 
 
 Under the hypotheses of  Theorem \ref{theo-3.2n}, for $k\ge 2$ 
  \begin{equation}
 \lim_{h\to 0}   \int_a^bH_{k}\(\frac{G(x+h)-G(x)}{\si (h)}
\)\,dx=0 
   \end{equation} 
whereas 
  \begin{equation}
 \lim_{h\to 0}   \int_a^b :\(\frac{G(x+h)-G(x)}{h}
\)^{k}:\,dx= \,\,:(G')^{k}:(I_{[a,b]}),
   \end{equation} 
   a well defined random variable, as we show in (\ref{3.7bs}) and Theorem \ref{theo-9.1}. Thus we see that when the hypotheses of  Theorem \ref{theo-3.2n} are satisfied it is quite natural to write the right-hand side of (\ref{mm5.22n}) in terms of Wick powers.
  
\medskip	 In Theorem \ref{theo-GCLT} we show that  for $\si^{2}$  relatively `large' at zero
\begin{equation}
    \int_a^bf\(\frac{G(x+h)-G(x)}{\si (h)}
\)\,dx-(b-a)Ef(\eta),
  \ee 
  divided  by it's variance, has a normal limit, in distribution, as $h\to 0$. We see this, heuristically, as a result of the fact that in these cases the increments of $G$ are only slightly correlated so that, writing the integral as a sum,  we are in the standard situation of a normal central limit theorem. (Note that when $\si^{2}$ is concave the increments of $G$ are negatively correlated.)
  
  On the other hand for $f$ and $\si^{2}(h)$ sufficiently smooth
  \begin{equation}
   \lim_{h\to 0}f\(\frac{G(x+h)-G(x)}{\si (h)}\)=f\({G'(x)\over \si'(0)}\)\qquad a.s.\label{6.10}
   \end{equation}
  as stochastic processes for $x\in [-T, T]$ for any $T>0$, where $G'$ is the actual derivative of $G$. In this case if we expand the right-hand side of (\ref{6.10}) in Hermite polynomials we get 
   \bea
&&\lim_{h\to 0} \int_a^bf\(\frac{G(x+h)-G(x)}{\si (h)}
\)\,dx \label{abstq}
 \\
&& \qquad= \sum_{j=0}^{\ff} (1/\si'(0))^{j}\,\,{E(H_{j}(\eta)
f(\eta))\over\sqrt {j!}}\,\int_{a}^{b}:(G'(x))^{j}: \, dx  \qquad a.s.\nn
\eea
  
 We now see that (\ref{mm5.22n}) lies somewhere between (\ref{j.1}) and 
 (\ref{abstq}). What distinguishes the hypotheses of Theorem \ref{theo-3.2n}   is that although $\si^{2}$ is not twice differentiable at zero, nevertheless
 \begin{equation}
   \int_{0}^{T}| (\si^{2})''(x)|^{j _{0}}\, dx<\ff.
   \end{equation} 
   We see in (\ref{mm5.22n}) what looks like the beginning of the power series expansion in (\ref{abstq}).  We see this even more dramatically in Example \ref{ex-5.2}, in which we show that for 
  $\si^{2}(u)\approx Cu^{2}\log^{2} 1/u$ and $ (\si^{2})''(u)\approx \log^{2} 1/u$,
  \be  \int_a^bf\(\frac{G(x+h)-G(x)}{\si (h)}
\)\,dx
\label{mm5.22nqq}  Ê \sim \sum_{j=0}^{\ff} (h/\si(h))^{j}\,\,{E(H_{j}(\eta)
f(\eta))\over\sqrt {j!}}\,\,:(G')^{j}:(I_{[a,b]})  
\ee
in $L^2$, as $h\to 0$.Ê
   
By considering a full range of Gaussian processes   we can appreciate how the asymptotic behavior of  (\ref{6.6}) changes as the increments variance of $G$ becomes smoother.

\medskip	  In Section \ref{sec-gd} we define the generalized derivative $G'$.  In
Section
\ref{sec-nonnor4} we construct the
$k$--th order Wick power process. This is used in Section \ref{sec7} to prove
Theorem
\ref{theo-3.2n} andÊ Corollaries \ref{theo-3.2}  and 
\ref{theo-3.1}.   Ê In Section \ref{sec-5} we give examples of Gaussian processes that satisfy the hypotheses of Theorem \ref{theo-3.2n}.  

\medskip	   There  are many papers about  non normal central limit theorems for nonlinear functionals of Gaussian processes. See for example  \cite[Dobrushin and Major]{DM},  \cite[Major]{Major},  
\cite[Taqqu]{Taqqu} and  \cite[Surgailis]{Surgailis}.   The focus of these papers differs significantly from what is considered in   this  paper. They consider long--range dependence and the limiting distributions that are obtained   are self--similar. In this paper we are concerned with local phenomena. The    generalized derivative $G'$ of the 
Gaussian process $G$,  appears  in the limit and it is clear from (\ref{3.7bq})  that the limiting distributions we obtain   are not, in general,  self--similar.  

Moreover because of the nature of the problems considered in the above references, they  only consider  weak convergence.  In contrast  we obtain   asymptotic expansions in $L^{2}$. This remark also applies to more recent results on the   non normal weak convergence of multiple Wiener--It\^{o} integrals; see, for example,   \cite[Nourdin and Peccati]{NP}, and the references therein.

\section{Generalized derivatives}\label{sec-gd}

The second condition in (\ref{mm5.17o})  implies that $G$ has a version with
continuous sample paths. (Clearly it implies that  $\si^2(h)\le Ch$, for
$h\in [0,h_0]$ for some constant $C$ and $h_0>0$. Therefore, continuity
follows from \cite[Lemma 6.4.6]{book}.) We work with this version.
  (It follows from the first condition in (\ref{mm5.17o}) that the paths of $G$ are not mean square   differentiable.)

 \begin{lemma}\label{lem-2.1}  Let  $G=\{G(x),x\geq 0\}$ be a mean zero  Gaussian  process
with stationary increments and $G(0)=0$, and with increments variance $\si ^2$ satisfying the second condition in (\ref{mm5.17o}).  If $\rho$ is locally integrable there exists a mean zero Gaussian field $\{G'(g), g\in\BB_{0}(R_{+})\}$    with covariance   
\begin{equation}
 E\(G'(g)G'(\wt g)\)=\int\int\rho (t-s)\,g( s)\,\wt g( t)   \,ds\,dt.\label{new.5}
\end{equation}
\end{lemma}
 
  We use the
following simple lemma which follows by simply doing the integration.

\bl \label{lem-3.1j}   Let $\phi$ be a symmetric function on $R^1$. Suppose
that $\phi''$ is locally integrable on $R^1$ and
$\phi(0)=\phi'(0)=0$. Then
\be
{1\over2}\int_{ a}^{ b}\int_{
a}^{b}     \phi''(x-y)   \,dx\,dy= \phi(b-a).\label{.m31qj}
\ee
\el 
%

\medskip	\noindent{\bf Proof of Lemma \ref{lem-2.1} }
It follows from Lemma \ref{lem-3.1j}  that
\be
 \si^2 (x)=\int_{ 0}^{ x}\int_{0}^{x}     \rho (t-s)   \,dt\,ds.\label{.m31qz}
\ee
Also, since $ G(0)=0$, $EG^{2}(x)=\si^2(x)$. Consequently for $x\leq y$
\begin{eqnarray}
&&
EG(x)G(y)\label{new.1}\\
&&\qquad    ={1\over2}\lc EG^{2}(x) +EG^{2}(y) 
-E\(G(x)-G(y)\)^{2}\rc\nonumber\\ &&\qquad    ={1\over2}\lc \si^2(x) +\si^2(y) 
-\si^2(y-x)\rc\nonumber\\ &&\qquad    ={1\over2}\int\int   \lc I_{\{[0,x]^{2}\}}
+I_{\{[0,y]^{2}\}}   -I_{\{[x,y]^{2}\}} \rc  \rho (t-s)   \,ds\,dt\nonumber\\
&&\qquad    ={1\over2} \lc \int_{0}^{x}\int_{0}^{y}    \rho (t-s)   \,ds\,dt 
+\int_{0}^{y}
\int_{0}^{x}    \rho (t-s)   \,ds\,dt    \rc\nonumber.
\end{eqnarray}
Since $\si^2$ is symmetric, so is $\rho$. Therefore  
\begin{equation}
EG(x)G(y)=\int_{0}^{x}\int_{0}^{y}    \rho (t-s)   \,ds\,dt.\label{new.2}
\end{equation}
 It follows from this that for $x'\leq x$, and $y'\leq y$
\begin{equation}
E\(G(x)-G(x')\)\(G(y)-G(y')\)=\int_{x'}^{x}\int_{y'}^{y}    \rho (t-s)   \,ds\,dt.\label{new.2star}
\end{equation}

Let $\mathcal{E}(R_+)$ be  the set of elementary  functions on $R_+$  of the form 
$g(x)=\sum_{i=1}^{n}g_{i}I_{\{(a_{i},b_{i}]\}}(x)$.
For such functions $g(x)$   we define the stochastic integral  
\begin{equation}
\int g(x)\,dG(x):=\sum_{i=1}^{n}g_{i}\(G(b_{i})-G(a_{i})\).\label{si.1}
\end{equation}
Note that by  (\ref{new.2star}), for  these functions,
\begin{eqnarray}
&&  \int\int\rho (t-s)\,g( s)\,g( t)   \,ds\,dt \label{si.2}\\
&&\qquad 
=\sum_{i,j=1}^{n}g_{i}g_{j}\int_{a_{i}}^{b_{i}}\int_{a_{j}}^{b_{j}}    \rho (t-s)   \,ds\,dt=E\(\int g(x)\,dG(x)\)^{2}\geq 0.\nn
\end{eqnarray}
It follows from this that the inner product
\begin{equation}
(g,\,\wt g)_{\mathcal{G}}:=\int\int\rho (t-s)\,g( s)\,\wt g( t)  
\,ds\,dt\label{fix.21}
\end{equation}
is positive definite on $\mathcal{E}(R_+)$.

Let $\mathcal{G}$ be the closure of
$\mathcal{E}(R_+)$   in the norm
\begin{equation}
\|g\|_{\mathcal{G}}=\(\int\int\rho (t-s)\,g( s)\,g( t)  
\,ds\,dt\)^{1/2}.\label{fix.21a}
\end{equation}
 Note that $\mathcal{G}$ is a Hilbert space. 
 It follows from (\ref{si.2}) that the stochastic integral extends from $\EE(R^{+})$ to a mean zero Gaussian field $\{G'(g), g\in\mathcal{G}\}$    with covariance   
\begin{equation}
 E\(G'(g)G'(\wt g)\)=(g,\,\wt g)_{\mathcal{G}}.\label{new.5j}
\end{equation}
  It is easy to see that $\mathcal{G}$ contains $\mathcal{B}_{0}(R_+)$.  \qed
  
  \br{\rm   There are several possible definitions of stochastic integrals for general Gaussian processes. See the discussion in \cite{SI} for the special case of fractional Brownian motion.
}
\er

We intend  the notation $G'$ to suggest the derivative.   If $G$ itself is
differentiable then $G'(g)$ could be written as  
\be
\int G'(x)g(x)\,dx,\label{92.1}
\ee
 in
which case the notation $G'(g)$ would be completely appropriate.  However,
even though the Gaussian processes that concern us are not differentiable
we may think of them as  having generalized derivatives for several reasons,
which  we give in the remainder of this section. 

\medskip

\bt \label{theo-si} Let $G$ be a Gaussian process of the type described in
Lemma \ref{lem-2.1}. Then  for any
$g\in
\mathcal{B}_{0}(R_+)$
\begin{equation} 
\lim_{h\rar 0}\int \({G( x+h)- G( x) \over h}\)\,g( x)\,dx=   G'(g)\hspace{.15 in}\mbox{ in }
L^{2}.\label{si.79}
\end{equation}
\et

\medskip
\Proof  
  Let
\begin{equation}
X_{ h}( g):=\int \({G( x+h)- G( x) \over h}\)\,g( x)\,dx. \label{nor4.4}
\end{equation}
We show that
\be
\lim_{h\to 0}E\(X_h(g)-G'(g)\)^2=0
\ee
by showing that all the terms of the expectation have the same limit as
$h\to0$.

Using the fact that $G( x+h)- G( x)=\int I_{\{(x,x+h]\}}(y)\,dG(y)$, it follows by
Fubini's Theorem
 and (\ref{new.5})    that
\begin{eqnarray}
&&E\lc X_h(g ) G'(g)\rc
\label{si.8}\\
&&\qquad    ={1 \over h}\int \,E\lc   \(G( x+h)- G( x)\)  G'(g) \rc \,g( x)\,dx\nonumber\\
&&\qquad    ={1 \over h}\int \,\lc  \int\int\rho (t-s)\,I_{\{(x,x+h]\}}(s)\, g( t)   \,ds\,dt  \rc \,g( x)\,dx\nonumber\\
&&\qquad    =  \int\int   \,\lc {1 \over h}\int_{s-h}^{s} \,g( x)\,dx\rc \rho (t-s)\,\, g( t)   \,ds\,dt.  \nonumber
\end{eqnarray}
By Lebesgue's theorem on differentiation
\be
\lim_{h\rar 0}{1 \over h}\int_{s-h}^{s}\,g( x)\,dx=g(s)\qquad \mbox{ for almost
all
$s$.}\label{Lebesgue}
\ee
Using this and the Dominated Convergence Theorem  we see that
\begin{equation}
\lim_{h\rar 0}E\(X_{ h}( g)G'(g)\)=\int\int\rho (t-s)\,g( s)\,g( t)   \,ds\,dt. \label{si.9}
\end{equation}

Considering (\ref{new.5}) we see that to complete the proof of this theorem
it suffices to show that
\begin{equation}
\lim_{h\rar 0}E\(X^{2}_{ h}( g)\)=\int\int\rho (t-s)\,g( s)\,g( t)   \,ds\,dt. \label{si.10}
\end{equation}
 Using (\ref{new.2star}) we have
\begin{eqnarray}
\lefteqn{ E\(X_{ h}( g)X_{ h'}(\wt g)\)
\label{new.4}}\\
&&= {1 \over h}{1 \over h'}\int\int E\(\(G(x+h)-G(x)\)\(G(y+h')-G(y)\)\)\,g( x)\,dx\,\wt g( y)\,dy
\nonumber\\
&&  = {1 \over h}{1 \over h'}\int\int\lc \int_{x}^{x+h}\int_{y}^{y+h'}    \rho (t-s)   \,ds\,dt\rc\,g( x)\,dx
\,\wt g( y)\,dy\nonumber\\
&& = \int\int\lc {1 \over h}\int_{t-h}^{t}\,g( x)\,dx\rc \lc {1 \over h'}\int_{s-h'}^{s} \,\wt g( y)\,dy \rc  \rho (t-s)   \,ds\,dt.\nonumber
\end{eqnarray}
It now follows from the  Dominated Convergence Theorem  
and (\ref{Lebesgue}),  that (\ref{si.10}) holds.  \qed

\br\label{almostsure}{\rm    When $g=I_{\{(a,b]\}} $, (\ref{si.79}) and the construction of $G'$ show that 
 \begin{equation}
\lim_{h\to 0}\int_{a}^{b} \({G( x+h)- G( x) \over h}\) \,dx =G'(I_{\{(a,b]\}})= G(b)-G(a)\hspace{.15 in}\mbox{ in }
L^{2}.\label{si.7a}
\end{equation}
It is easy to see that this limit actually holds almost surely.  Since $G$ has continuous paths almost surely,
\begin{eqnarray}
   &&  \lim_{h\to 0}\int_{a}^{b} \({G( x+h)- G( x) \over h}\) \,dx  \label{snor}\\
  &&\qquad=  \lim_{h\to 0} {1 \over h}\lc \int_{a}^{b} G( x+h)\,dx -\int_{a}^{b} G( x)\,dx\rc \nn\\
  &&\qquad = \lim_{h\to 0} {1 \over h}\int_{b}^{b+h} G( x)\,dx -{1 \over h}\int_{a}^{a+h} G( x)\,dx = G(b) -G(a)\qquad
\mbox{a.s.}\nn
\end{eqnarray}
More generally,  for all $g\in \mathcal{E}(R_+)$ we actually have almost sure
convergence in (\ref{si.79}).
}
\er

Finally we note that we can consider $G$ to be a   (random) distribution
defined by
\begin{equation}
G(f)=\int G(x)f(x)\,dx, \hspace{.2 in}  f\in C_{0}^{\ff}(R_+)\label{d.1}.
\end{equation}
In this case $G$ has a distributional derivative $DG$.  Using the fact that
 for any   $f\in C_{0}^{\ff}(R_+)$, $(f(x+h)-f(x))/
 h$ converges to $f'(x)$ uniformly, we
have
\begin{eqnarray}
DG(f)&:=&-G(f') 
\label{d.2}\\
&=&\lim_{h\rar 0}\int G(x){f(x-h)-f(x) \over h}\,dx \nonumber\\
&=&\lim_{h\rar 0}\int {G(x+h)-G(x) \over h}f(x)\,dx\nonumber
\end{eqnarray}
almost surely. 
Therefore, it follows from  Theorem \ref{theo-si}  that for $f\in
C_{0}^{\ff}(R_+)$, 
$G'(f)\stackrel{L^{2}}{=}DG(f)$.

\section{Wick powers of generalized derivatives}\label{sec-nonnor4}

 Let   $(X,Y)$ be a
two dimensional Gaussian random variable.     
By   \cite[Theorem 3.9]{Janson}  
\begin{equation} E(:X^{ k}::Y^{  j}:)=k!(E( XY))^{ k}\de_{
k,j}.\label{18.3a}
\end{equation}
($:X^{ k}:$ is defined in (\ref{defin}).)      It follows from   (\ref{0.9}) and (\ref{18.3a}) that if  $X$ and $Y$  are $N(0,1)$ and $(X,Y)$ is a
two dimensional Gaussian random variable then
\begin{equation} E(H_{ k}( X)H_{ j }( Y))=(E( XY))^{  k}\de_{
k,j}.\label{18.3}
\end{equation}
 
  We   say that a function $\varrho( x)$ is weakly positive definite if  \begin{equation}
\int\int \varrho( s-t)g( s)g( t)\,ds\,dt\geq 0 
\end{equation}
for all $g\in \mathcal{B}_{0}(R_+)$.
Let  $\varrho( x)$ be a 
  symmetric, weakly positive definite function that is locally integrable on 
  $R^{1}$.    Consider the    mean zero Gaussian field   
$\FF=\{\FF(g), g\in \mathcal{B}_{0}(R_+) \}$ 
   with
 covariance  
\begin{equation}
E(\FF( f)\FF( g))=\int\int \varrho( s-t)f( s)g( t)\,ds\,dt\qquad
f,g\in\mathcal{B}_{0}(R_+)\label{nor4.1}.
\end{equation}
(We are particularly interested in the case in which
  $\varrho( 0)=\ff$,  in which case it is not the covariance of a
stationary Gaussian process.)

 Let $f_{\de}(s)$ be a continuous positive symmetric function on
$(s,\de)\in R_+ \times (0,1]$, with support in the ball of radius $\de$ centered at the origin, with
$\int f_{\de}(y)\,dy=1$. That is,
$f_{\de}$ is a   continuous  approximate identity.  Set
$f_{x,\de}(s)=f_{\de}(s-x)$. 
 
 Assume that 
 \be
\varrho\in L^{ k}_{loc}(R^{1}).\label{3.4}
 \ee
  We now define,  what
we call, the
$ k$-th Wick power  Gaussian chaos associated with $\FF$.  

\begin{lemma}\label{lem-3.1} Let $\{f_{\de},\de\in(0,\de_{0}]\}$ be a family of approximate identities and assume that $\varrho$ is a 
   symmetric, weakly positive definite function that satisfies (\ref{3.4}). Then for all $g\in\BB_{0}(R_{+})$
\be
:\FF^k:(g)\stackrel{def}=\lim_{\de\rar 0}\int
:(\FF ( f_{x,\de}))^k: g(x)\,dx\qquad \mbox{exists in $L^2$}.\label{wwc}
\ee
and
\be
E( :\FF^k:(g) )^2=k!\int\!\!\int \varrho^{k}(x-y) g(x)g(y)\,dx\,dy\label{3.7b}.
\ee
 
  \end{lemma}

\Proof Consider the
mean zero Gaussian process
$\{\FF( f_{x,\de})\,, (x,\de)\in R_+ \times (0,1]\}$ with covariance
 \bea
E(\FF( f_{x,\de})\FF( f_{y,\de'}) ) &=&\int\int
\varrho( x'-y')f_{x,\de}(x') f_{y,\de'}(y')\,dx'\,dy'\label{3.5}\\ &=&\int\int
\varrho( x'+x-y'-y)f_{\de}(x') f_{\de'}(y')\,dx'\,dy' \nn\\
&\stackrel{def}{=}&\varrho_{\de,\de'}(x,y).\nn \eea 
It follows from (\ref{18.3a})  that  
\be E(:(\FF( f_{x,\de}))^{k}:\,:(\FF( f_{y,\de'}))^{k}:)=k!(\varrho_{\de,\de'}(x,y))^{k}.
\label{11.2}
\ee

  Let $g\in \BB_{0}(R_{+}) $.  It follows
from  (\ref{3.5}), (\ref{11.2}) and Fubini's theorem that
\bea \lefteqn{ E\(\int \!\!\int :(\FF ( f_{x,\de}))^k:\,:(\FF ( f_{y,\de'}))^k:
g(x)g(y)\,dx\,dy\)
\label{3.7a} } \\ &&
=k!\int\!\!\int (\varrho_{\de,\de'}(x,y))^{k}g(x)g(y)\,dx\,dy\nn 
\\  &&
=k!\int\!\!\int\(\int\dots\int\) \prod_{j=1}^{k}\varrho (x+v_j-y-w_j) \prod_{j=1}^{k}
f_{\de}(v_j) f_{\de'}(w_j)\,dv_j\,dw_j \nn\\ 
 &&\hspace{3.5in} g(x)g(y)\,dx\,dy\nn\\&&
=k!\int\dots\int\(\int\!\!\int \prod_{j=1}^{k}
\varrho(x-y+v_j-w_j) g(x)g(y)\,dx\,dy\)\nn\\
&&\hspace{3in} \prod_{j=1}^{k}
f_{\de}(v_j) f_{\de'}(w_j)\,dv_j\,dw_j.\nn
 \eea 
Since $\varrho\in L^{k}_{loc}(R^{1})$ and translation is continuous in $
L^{k}_{loc}(R^{1})$, the double integral in parentheses immediately above is
continuous in
$(v_1-w_1,\ldots, v_{k}-w_{k})$ and goes to 
\be
\int\!\!\int\varrho^{k}(x-y) g(x)g(y)\,dx\,dy\label{3.7aa}
\ee as $\sup_{1\le j\le n}|v_j-w_j|\to 0$.  
 Consequently
\bea
&&\lim_{\de,\de'\to0}E\(\int \!\!\int:(\FF ( f_{x,\de}))^k:\,:(\FF (
f_{y,\de'}))^k: g(x)g(y)\,dx\,dy\)\label{92.5}\\
&&\qquad \quad=k!\int\!\!\int\varrho^{k}(x-y) g(x)g(y)\,dx\,dy\nn.
\eea
  It follows from this that 
\be 
\lim_{\de,\de'\to0}E\( \int:(\FF ( f_{x,\de}))^k:\,  g(x) \,dx - \int:(\FF (
f_{x,\de'}))^k:\,  g(x) \,dx\) ^2=0.
\ee 
  This implies (\ref{wwc}). The relation in (\ref{3.7b}) follows from (\ref{92.5}).\qed  



 \begin{remark} {\rm 
 Suppose that
$\wt{\FF}$ is a   mean zero Gaussian field, with covariance 
$\wt\rho \in L^{k}_{loc}(R^{1})$, and that
 $\FF$ and $\wt{\FF}$  are jointly Gaussian with
\begin{equation}
E(\FF(f)\wt{\FF}(f'))=\int\int \psi( x-y)f( x)f'( y)\,dx\,dy\label{.mm1.1}
\end{equation}
for some $\psi\in L^{k}_{loc}(R^{1})$.   If we return to (\ref{3.5}) and replace $\FF( f_{y,\de'}))$ by $\wt \FF( f_{y,\de'}))$ and continue the argument in the proof of  Lemma \ref{lem-3.1}   we see
that 
 \be
E (:\FF^{k}:(g):\wt\FF^{k}:(g')) =k!\int\!\!\int
 \psi^k( x-y) g( x)g'( y)\,dx\,dy.\label{3.11a} 
\ee 
 }\end{remark}

 \br\label{rem-3.2}{\rm    Although we say that we are particularly interested in the case in which
  $\varrho( 0)=\ff$ in (\ref{nor4.1}), Lemma 3.1 also applies when $\varrho$ is the covariance  of a stationary Gaussian process. Given 
 a mean zero  stationary Gaussian process  $\wt
G=\{\wt G(x),x\ge 0 \}$, with  continuous covariance
$\varphi(s)$,  we can define  a Gaussian field $\GG=\{\GG(g),g\in \BB_0(R_+ )\}$ by  
\begin{equation}
\GG(f)=\int \wt G(x)\,f(x)\,dx\label{3.11b}.
\end{equation}
Clearly
\begin{equation}
E(\GG( f)\GG( g))=\int\int \varphi( s-t)f( s)g( t)\,ds\,dt\qquad
f,g\in\mathcal{B}_{0}(R_+)\label{nor4.1h}.
\end{equation}
  It follows from Lemma \ref{lem-3.1} that    we can construct a $k$--th
order Wick power chaos  $\GG^k=\{ :\GG^k:(g), g\in \BB_0(R_+ )\}$ with 
\be
E( :\GG^k:(g) )^2=k!\int\!\!\int \varphi^{k}(x-y)
g(x)g(y)\,dx\,dy\label{3.7bqas}.
\ee 
However, we do not really need Lemma \ref{lem-3.1} when we are dealing with a mean zero stationary Gaussian process, since we can simply form the $k$--th
order Wick power chaos $\wt\GG^{k}=\{\wt\GG^{k}(g),g\in \BB_0(R_+ )\}$ by setting
\begin{equation}
:\wt \GG^{k} :(g)=\int :\(\wt G(x)\)^{k}:\,g(x)\,dx.\label{3.11ca}
\end{equation}

It  is    easy to see that these two processes, $ \GG^{k}$  and $\wt\GG^{k}$, are equivalent, (in $L^{2}$).   To do this we now show that 
\be
\lim_{\de\rar 0}\int
:(\GG ( f_{x,\de}))^k: g(x)\,dx=\int :\(\wt
G(x)\)^{k}:\,g(x)\,dx,\label{3.11caq}\qquad \mbox{ in $L^2$}.
\ee

  Note that by (\ref{3.11b}) and the fact that $\wt G$ has covariance $\varphi$,
 \bea
E(\GG( f_{x,\de}) \wt G(y)) &=&\int\int
\varphi ( x'-y)f_{x,\de}(x')\,dx'\label{3.5f}\\ &=&\int\int\varphi ( x'+x-y)f_{\de}(x') \,dx' \nn\\
&\stackrel{def}{=}&\varphi _{\de}(x,y).\nn \eea 
Therefore, it follows from (\ref{18.3a})  that  
\be E(:(\GG( f_{x,\de}))^{k}:\,:(\wt G(y))^{k}:)=k!(\varphi_{\de}(x,y))^{k}.
\label{11.2f}
\ee
We use this in place of (\ref{11.2}) and continue with the  
 argument in (\ref{3.7a})--(\ref{92.5}), with obvious modifications, to see that 
\bea
&&\lim_{\de\to0}E\(\int \!\!\int:(\GG ( f_{x,\de}))^k:\,:(\wt G(y))^k:
g(x)g(y)\,dx\,dy\)\\
&&\qquad \quad=\int\!\!\int\varphi ^{k}(x-y) g(x)g(y)\,dx\,dy\nn.
\eea
Using this, (\ref{92.5}) with $\varrho$ replaced by $\varphi$,  and the obvious fact that 
\bea 
&& E\(\int \!\!\int :\(\wt G(x)\)^{k}\, :\(\wt G(y)\)^{k}:
g(x)g(y)\,dx\,dy\) \\
&&\qquad  =k!\int\!\!\int\varphi ^{k}(x-y) g(x)g(y)\,dx\,dy ,\nn
\eea 
we get (\ref{3.11caq}).  }\end{remark}

  \begin{remark}\label{rem-3.3}{\rm In this paper we consider Wick powers of Gaussian fields,  $\{:(G')^{k}:(g), g\in \BB_{0}(R_{+})\}$
It is well known that $:(G')^{k}:(g)$ can also be expressed as a multiple   Wiener-It\^{o} integral. ( See, e.g.  \cite{Major}.) We briefly explain this  for the benefit of those familiar with multiple   Wiener-It\^{o} integrals: 

Since $\rho (x)$ is symmetric and weakly positive definite, it follows from the Bochner--Schwartz Theorem that  $\rho (x)=\int e^{i\la x}\,d\mu(\la)$ for some positive   Radon measure $\mu$.  When $\rho (0)=\ff$,    $\mu$ is not a finite measure.

Let  $Z_{\mu}$ be the (complex valued) Gaussian random spectral measure corresponding to $\mu$.
Then 
\begin{equation}
:(G')^{k}:(g)=\int\cdots\int \wh g(\la_{1}+\cdots+\la_{k})\,dZ_{\mu}(\la_{1})\cdots\,dZ_{\mu}(\la_{k})\label{dob1}
\end{equation}
where $\wh g$ is the Fourier transform of $g$.
(This is the end of Remark \ref{rem-3.3}.)
 }\end{remark}

  \medskip  We now apply the above results about constructing   Gaussian
chaoses to the processes that concern us.  In Lemma \ref{lem-2.1} we
define  the Gaussian field
$\{G'(g), g\in\BB_0(R_+ )\}$.  
When
$\rho\in L^{k}_{loc}(R^{1})$ the procedure that leads to (\ref{wwc}) and
(\ref{3.7b}) enables us to define 
$k$-th Wick power chaos
 \be
:(G')^k:(g)=\lim_{\de\rar 0}\int
:(G' ( f_{x,\de}))^k: g(x)\,dx\label{wwcs}
\ee
  as a limit in $L^2$, with
\be
E( :(G')^k:(g) )^2=k!\int\!\!\int \rho^{k}(x-y) g(x)g(y)\,dx\,dy\label{3.7bs}.
\ee
Note also that the Gaussian field  $X_{ h}( g)$ defined in (\ref{nor4.4}) is of
the form of (\ref{3.11b}).
 Therefore, it follows from (\ref{3.11ca}) that
\begin{equation}
:X_{h}^{  k}:( g)=\int
:\({G( x+h)- G( x) \over h}\)^{ k}:\,g( x)\,dx.\label{nor4.4m}
\end{equation}

The next theorem is a critical result in this paper.  

\bt \label{theo-9.1}  LetÊ
$G=\{G(x),x\in R_+\}$, $G(0)=0$, be a mean zero Gaussian processÊ with
stationary increments. Let $\rho$ Êbe as defined in (\ref{rho}) and assume that   $\rho^{ k }( x)$ is  locally integrable and that   $\rho(|x|)$ is
bounded  on $[\de,M]$ for each $0<\de<M<\ff$ . Then for all $g\in\BB_0(R_+ )$,
\begin{equation}
  \lim_{h\to 0} \int
:\({G( x+h)- G( x) \over h}\)^{ k}:\,g( x)\,dx= \,\, :(G')^{ k} :( g) \qquad \label{nor4.6}
\mbox{in
$\quad L^2$}.
\ee
 \et

\Proof By (\ref{new.4})  
\begin{eqnarray}\lefteqn{
E (X_{ h}( g)X_{ h'}(  g))
\label{nor4.3c}}\\
&& = \int \int \({1 \over h}\int_{ x}^{ x+h}
{1 \over h'}\int_{ y}^{ y+h'}\rho( s-t)\,dt \,ds\)\,g( x)g( y)\,dx\,dy. \nonumber
\end{eqnarray}
Consequently by (\ref{.mm1.1}) and (\ref{3.11a})
\begin{eqnarray}\lefteqn{
E (:X_{ h}^{k}:( g):X^{k}_{ h'}:(  g))
\label{nor4.3caa}}\\
&& = k!\int \int \({1 \over h}\int_{ x}^{ x+h}
{1 \over h'}\int_{ y}^{ y+h'}\rho( s-t)\,dt \,ds\)^{k}\,g( x) g( y)\,dx\,dy.
\nonumber
\end{eqnarray}

  In a similar vein by  (\ref{si.79}) and (\ref{new.4}),  
Lebesgue's Theorem and a change of variables 
\begin{eqnarray}
\lefteqn{
E (X_{ h}( g)G'( g))
\label{nor4.3b}}\\
&& = \lim_{h'\rar 0} E (X_{ h}( g)X_{ h'}( g))\nonumber\\
&& =  \lim_{h'\rar 0} \int\int\lc {1 \over h}\int_{t-h}^{t}\,g( x)\,dx\rc \lc {1 \over h'}\int_{s-h'}^{s} \,g( y)\,dy \rc  \rho (t-s)   \,ds\,dt\nn\\
&& = \int\int\lc {1 \over h}\int_{t-h}^{t}\,g( x)\,dx\rc   \rho (t-s)   \, g( s)\,ds\,dt \nonumber\\
&& = \int \int \({1 \over h}\int_{ x}^{ x+h}\rho( s-y) \,ds\)\,g( x) g( y)\,dx\,dy. \nonumber
\end{eqnarray}
Therefore, by (\ref{.mm1.1}) and (\ref{3.11a}) 
\bea
&&E (:X^{k}_{ h}:( g):(G')^{k}:(  g))\label{.mm1.2}\\
&&\qquad =k! \int \int \({1 \over h}\int_{ x}^{
x+h}\rho( s-y) \,ds\)^{k}\,g( x)g (y)\,dx\,dy\nn\\
&&\qquad =A_{h,\de}+B_{h,\de}
\nn
\eea 
where 
\begin{equation}
A_{h,\de}=k! \int \int_{|x-y|<\de} \({1 \over h}\int_{ x}^{
x+h}\rho( s-y) \,ds\)^{k}\,g( x) g (y)\,dx\,dy\label{til.1}
\end{equation}
and
\begin{equation}
B_{h,\de}=k! \int \int_{|x-y|\geq \de} \({1 \over h}\int_{ x}^{
x+h}\rho( s-y) \,ds\)^{k}\,g( x) g (y)\,dx\,dy.\label{til.2}
\end{equation}
Fix $\de>0$. Using the fact that $\rho$ is  bounded away from
the origin,   the Dominated Convergence Theorem, and Lebesgue's theorem
on differentiation,  we see that
\begin{equation}
\lim_{h\rar 0}B_{h,\de}=k! \int \int_{|x-y|\geq \de} \rho^{k}( x-y) \,g( x)g (y)\,dx\,dy.\label{til.3}
\end{equation}

On the other hand, using  
the H\"older or Jensen inequality, we see that  for $h\leq \de$
\begin{eqnarray} 
|A_{h,\de}|&&\le k!
\int\!\!\int_{|x-y|<\de} \({1 \over h}\int_{ x}^{ x+h}|\rho^{k}( s-y)| \,ds\)|g(x)|\,| g(
y)|\,dx\,dy\nn\\ &&\leq  k!\int\!\!\int_{|s-y|<2\de} |\rho^{k}( s-y)|\({1 \over h}\int_{s-h}^{
s}|g(x)|\,dx\) |g( y)|\,dy \,ds. \nn\\ &&\leq  C\int_{|s |<2\de} |\rho^{k}( s )|  \,ds. \label{nor4.9}
\end{eqnarray}
Since by assumption $\rho^{k}( s ) $ is locally integrable  we can make this arbitrarily small by choosing $\de>0$ sufficiently small. Thus we have shown that
\be
\lim_{h\to0}E( :X_{h}^{ k}:( g)\, :(G') ^{ k}:( g))=  k! \int\!\!\int \rho^{k}( s-y)
g(s)  g( y)\,dy \,ds.\label{.3.1q}
\ee

 Similar reasoning shows that   
 \be
\lim_{h\to0}E( :X_{h}^{ k}:( g) :)^2 =   k!\int\!\!\int \rho^{k}(
s-y) g(s)  g( y)\,dy \,ds\label{.3.1qq}
\ee
 for all   $g\in\BB_{0}(R_{+})$.    Using (\ref{.3.1q}),  (\ref{.3.1qq}) and
(\ref{3.7bs}) we get (\ref{nor4.6}).
 \qed

 \br{\rm 
In Section \ref{sec-gd} we explain why we think of the field $G' $ as a
generalized derivative of the Gaussian 
 process $G=\{G(x),x\in R^1\}$.   In (\ref{wwcs})  we construct  the $k$-th
Wick power chaos $:(G')^k:(g)$.
  When
$G$ itself is mean square differentiable, i.e. when
\be
\lim_{h\to 0}E \({G( x+h)- G( x) \over h}\)^2=\lim_{h\to 0}{\si^2(h)\over
h^2}:=\frac{1}{2}(\si^2)''(0)<\ff,
\ee
$\{{d\over dx}G(x), x\in R^1\}$ is a stationary Gaussian process with covariance
$\frac{1}{2}(\si^2)''(x-y) $. In this case, as we show  in (\ref{3.11ca}),
\begin{equation}
:(G')^k:(g)=\int :\({d\over dx}G(x)\)^{k}:\,g(x)\,dx.\label{3.11cad}
\end{equation}
This further motivates the description of $G' $ as a
generalized derivative.

  However, for Gaussian processes satisfying
(\ref{mm5.17o}),
$\lim_{h\to 0}{\si^2(h)\over
h^2}=\ff$,  and consequently  
${d\over dx}G(x)
$ is not a stochastic process.    (In these cases, formally taking $g$
in (\ref{new.5}) to be the delta `function'
$\de_{x}$,  gives    $E(G'(x))^{2}=\rho(0):=\lim_{h\to 0}{\si^2(h)\over
h^2}=\ff$.)

}
\er

\section{$L^{2}$ asymptotic expansion}\label{sec7}

For each
$h$ we consider the symmetric positive definite kernel
\bea
\tau_{h}(x,y)&=&{ 1\over \sigma^{ 2}(h)}E(G(x+h) - G(x))(G(y+h) -
G(y))
\label{ag2.3}\\
      &=&{ 1\over 2\sigma^{2}(h)}\(\si^2(x-y+h)+ \si^2(x-y-h)-2\si^2(x-y)
  \)\nn\\ &:=&\tau_h(x-y)=\tau_h(y-x),\nn
\eea  
 (see (\ref{new.1})). Note that since $G$ has stationary
increments it follows from the  Cauchy--Schwarz inequality that
\begin{equation} |\tau_{h}(x-y)|\leq 1\qquad \hspace{
.2in}\forall\,x,y\in R^{ 1}.\label{ag2.3t}
\end{equation}

\medskip  To continue we need some estimates of the integrals of powers of $\tau_{h}$.

\bl \label{lem-3.2}Suppose that  $\si^2$
satisfies (\ref{rv})--(\ref{rv2}) and  
   $ \rho(s)$ is locally integrable   and is bounded in compact neighborhoods excluding  the origin. 
Then
\be
\lim_{h\to 0} {\int_{ a}^{ b}\!\!\int_{ a}^{ b}
 \tau _{h}(x-y) 
\,dx\,dy\over  h^2/ \si^2(h) }=\si^2(b-a) \label{mm5.2}
\ee
and
\be
\lim_{h\to 0} {\int_{ a}^{ b}\!\!\int_{ a}^{ b}
 \tau^2_{h}(x-y) 
\,dx\,dy\over  \int_{ a}^{ b}\!\!\int_{ a}^{ b}
 \tau_{h}(x-y) 
\,dx\,dy}=0.\label{mm5.2h}
\ee

More generally if,   in addition,  $\rho^k(s) $ is locally
integrable for some integer $k\ge 1$ and  
\be
h=o\( h^2/\si^2(h) \)^{k}\label{regularityb}
\ee
then
\be
\lim_{h\to 0} {\int_{ a}^{ b}\!\!\int_{ a}^{ b}
 \tau^k_{h}(x-y) 
\,dx\,dy\over  (h^2/ \si^2(h))^k }=\int_{ a}^{ b}\!\!\int_{ a}^{
b}\rho^k(x-y)\,dx\,dy \label{mm5.2q}
\ee
and
\be
\lim_{h\to 0} {\int_{ a}^{ b}\!\!\int_{ a}^{ b}
 \tau^{k+1}_{h}(x-y) 
\,dx\,dy\over  \int_{ a}^{ b}\!\!\int_{ a}^{ b}
 \tau^k_{h}(x-y) 
\,dx\,dy}=0.\label{mm5.2hw}
\ee
\el

\Proof   It follows from  
(\ref{ag2.3})  and  (\ref{nor4.4}) that 
\bea
&&{\int_{ a}^{ b}\!\!\int_{ a}^{ b}
 \tau_{h}(x-y) 
\,dx\,dy\over  h^2/ \si^2(h) }\label{ag2.3v}\\
&&\qquad =   \int_{ a}^{ b}\!\!\int_{ a}^{ b}
{\si^2(x-y+h)+ \si^2(x-y-h)-2\si^2(x-y)
\over 2h^2}
\,dx\,dy\nn \\
&&\qquad   =   E(X _{h}(I_{[a,b]}))^{2} .\nn
\eea
By  (\ref{si.7a}), 
\be 
 \lim_{h\to 0} E(  X_{h} (I_{[a,b]})  )^2 =E(G'(I_{[a,b]}))^2=E(G (b
)-G(a))^2 =\si^2(b-a).
\ee 
Thus we get (\ref{mm5.2}).

The statement in (\ref{mm5.2q}) follows as above using   (\ref{nor4.3caa}) and then (\ref{nor4.6}) which
implies that
\begin{equation}
  \lim_{h\to 0} E( :X_{h}^{ k}:(I_{[a,b]}) :)^2\,= E\(:(G')^{ k}:( I_{[a,b]}) \)^2
.
\end{equation}

We now obtain (\ref{mm5.2hw}), which,  considering (\ref{mm5.17o}), includes  (\ref{mm5.2h}). We are given that
$\rho^k(s)$ is locally integrable. Suppose that
$\rho^{k+1}(s)$ is also locally integrable.   Then, by  (\ref{mm5.2q}) 
\be
\lim_{h\to0}\({\si^2(h)\over h^2}\)^{k+1}\int_{ a}^{ b}\!\!\int_{ a}^{ b}
 \tau^{k+1}_{h}(x-y)  
\,dx\,dy = 
 \int_{ a}^{ b}\!\!\int_{ a}^{
b}\rho^{k+1}(x-y)\,dx\,dy.
\ee
 The statement in (\ref{mm5.2hw}) clearly follows from this, (\ref{mm5.17o})
and (\ref{mm5.2q}).

  Suppose  $\rho^{k+1}(s)$ is not   
integrable over neighborhoods of the origin.  By a change of variables, and with $c=b-a$,
\bea
&&\int_{ a}^{ b}\!\!\int_{ a}^{ b}
 (\tau_{h}(x-y) )^{k+1}
\,dx\,dy\label{fin2} \\
&&\qquad = 2  \int_0^c 
\({ \si^2(s+h)+ \si^2(s-h)-2\si^2(s)\over 2\si^2(h)}\)^{k+1}
 (c-s)\,ds\nn \\
&&\qquad \le 16ch+   \int_{8h}^c \(
{ \si^2(s+h)+ \si^2(s-h)-2\si^2(s)\over 2\si^2(h)}\)^{k+1}
 (c-s)\,ds\nn,
\eea 
where, for the last line we use (\ref{ag2.3t}). Also
\bea
&&\int_{8h}^c 
\({ \si^2(s+h)+ \si^2(s-h)-2\si^2(s)\over \si^2(h)}\)^{k+1}
 (c-s)\,ds\label{final}\\
&&\qquad =\({h^2\over \si^2(h)}\)^{k+1}\int_{8h}^c 
\({ \si^2(s+h)+ \si^2(s-h)-2\si^2(s)\over h^2}\)^{k+1}
 (c-s)\,ds\nn\\
&&\qquad \le c \({h^2\over \si^2(h)}\)^{k+1}\int_{8h}^c 
\({ \si^2(s) \over s^2}\)^{k+1}
 \,ds \nn,
\eea 
 where, for the last line we use (\ref{mm5.9}).  
Let $a>0$. Using (\ref{mm5.9}) again we see that
\bea
&&\int_{a}^c  \rho^{k+1}(s)\,ds\\
&&\qquad =\frac{1}{2^{k+1}}\int_{a}^c \lim_{h\to0}
\({ \si^2(s+h)+ \si^2(s-h)-2\si^2(s)\over h^2}\)^{k+1}
  \,ds\nn\\
&&\qquad \le C\int_{a}^c\({\si^2(s)\over
s^2}\)^{k+1}\,ds.\nn
\eea
Consequently, since $\rho^{k+1}(s)$ is not  locally
integrable,   the final integral in (\ref{final}) goes to infinity as
$h\downarrow0$. Since $\si^2(s)$ is regularly varying at zero, this means that 
$\({ \si^2(s) \over s^2}\)^{k+1}$ is regularly varying at zero with index
less than or equal to $-1$.

 Suppose that its index is equal to $-1$.   This implies that $\({h^2\over
\si^2(h)}\)^{k}$ is regularly varying with index equal to $k/(k+1)<1$ and that
the integral in the last line of (\ref{final}) is slowly varying. Consequently
\be
\({h^2\over \si^2(h)}\)^{k+1}\int_{8h}^c 
\({ \si^2(s) \over s^2}\)^{k+1}\label{fin3}
 \,ds 
\ee
is regularly varying with index equal to 1.  Taking
(\ref{mm5.2q}), (\ref{fin2}) and (\ref{final}) into account we see that
(\ref{mm5.2hw}) holds in this case. 

Finally, suppose that the index of $\({ \si^2(s) \over s^2}\)^{k+1}$ is less than
$-1$. In this case 
\be
\({h^2\over \si^2(h)}\)^{k+1}\int_{8h}^c 
\({ \si^2(s) \over s^2}\)^{k+1}\label{fin4}
 \,ds \sim Ch
\ee
at 0, for some constant $C$.  Taking 
(\ref{mm5.2q})  and (\ref{regularityb}) into account we again
get (\ref{mm5.2hw}). \qed

\medskip	\noindent{\bf Proof  of Theorem \ref{theo-3.2n} }  It follows from 
(\ref{a18.1}) and (\ref{am.m1}) that 
\bea
&&Ê \int_a^b f \(\frac{G(x+h)-G(x)}{\si (h)}
\)\,dx\label{mm5.19j}\\
&&\qquad =\sum_{j=0}^{j_0}E(H_{j}(\eta) f(\eta))
\int_a^bH_{j}\(\frac{G(x+h)-G(x)}{\si(h)}
\)\,dx\nn\\
&&\qquad
 \qquad +\sum_{ j=j_{0}+1}^{ \ff}a_{
j}\int_a^b H_{j}\(\frac{G(x+h)-G(x)}{\si (h)}
\)\,dx\nn.
\eea
Denote the last line in (\ref{mm5.19j}) by $\ZZ(h)$.
Using (\ref{18.3}) and then (\ref{ag2.3t}) { we have
\bea
E\ZZ^2(h)&=& \sum_{ j=j_{0}+1}^{ \ff}a_{
j}^2 \int_{ a}^{ b}\!\!\int_{ a}^{ b}
 \tau^j_{h}(x-y) 
\,dx\,dy\label{nach}\\
&\le&  \int_{ a}^{ b}\!\!\int_{ a}^{ b}
 \tau^{j_{0}+1}_{h}(x-y)  
\,dx\,dy    \sum_{ j=j_{0}+1}^{ \ff}a_{
j}^2.\nn
\eea
It follows easily from (\ref{nor14.5a}) with $j_0\ze<1$  and 
(\ref{.m31qz}) that  (\ref{regularityb}) holds with $k$ replaced by $j_0$.
It then follows from (\ref{mm5.2hw}) and (\ref{mm5.2q}) that  the last line in (\ref{nach}) is $o\({h\over\si (h)}\)^{j_0}$.   
Since
\begin{eqnarray}
&&\sum_{j=0}^{j_0}E(H_{j}(\eta) f(\eta))
\int_a^bH_{j}\(\frac{G(x+h)-G(x)}{\si(h)}
\)\,dx
\label{asym.1}\\
&&\qquad ÊÊ =\sum_{j=0}^{j_0}(h/\si(h))^{j}\,\,{E(H_{j}(\eta)
f(\eta))\over\sqrt {j!}}\,\,
\int_a^b:\(\frac{G(x+h)-G(x)}{h}
\)^{j}:\,dx,\nonumber
\end{eqnarray}
 we see that (\ref{mm5.22n}) follows from (\ref{.06}) in the next lemma.  \qed

 \bl \label{calc}Let  
$G=\{G(x),x\in R_+\}$, $G(0)=0$, be a mean zero Gaussian process  with
stationary increments    that satisfies the hypotheses of Theorem  \ref{theo-3.2n}.
Then  for $1\le j\le j_0$   and   $g\in \BB_{0}(R_{+})$ 
\be
\| :X_{h}^{j}:(g)-:(G') ^{ j}:( g)\|_{2}\le C(|h|\varphi^j(h))^{1/2}\label{.03}
\ee
and  
\be
\lim_{h\to 0}\| :X_{h}^{j}:(g)-:(G') ^{ j}:( g)\|_{2} \({\si^2(h)\over
h^2}\)^{(j_0-j)/2}=0.\label{.06}
\ee

\el

\Proof 
To obtain (\ref{.03})
 we use (\ref{nor4.3caa}) for $\|:X_{h}^{j}:(g) \|_{2}^{2}$, (\ref{3.7b}) for 
$\|:(G') ^{ j}:( g)\|_{2}^{2}$ and  (\ref{.mm1.2}) for $E\(X_{h}^{j}:(g):(G') ^{ j}:( g)\)$  
to see that  
 \begin{eqnarray}
\lefteqn{\frac{1}{j!}\|:X_{h}^{j}:(g)-:(G') ^{ j}:( g)\|_{2}^{2}
\label{nor14.7bk}}\\
&&Ê =\int\int \lc
\({1 \over h}\int_{ x}^{ x+h}{1 \over h}
\int_{ y}^{ y+h}\rho( s-t)\,dt \,ds\)^{j} - \({1 \over h}\int_{ x}^{
x+h}\rho( s-y)
\,ds\)^{j}\right.\nn\\ &&\quad\quad\quadÊÊ - \left.
\(Ê {1 \over h}\int_{ y}^{
y+h}\rho(x-t)
\,dt\)^{j}
+\(\rho( x-y) \)^{j}\rc g(x) g( y)\,dx\,dy. \nonumber
\end{eqnarray}
Set $z=x-y$. We write 
\begin{eqnarray}
&&{1 \over h}\int_{ x}^{ x+h}{1 \over h}
\int_{ y}^{ y+h}\rho( s-t)\,dt \,ds
\label{ct.1}\\
&&\qquad   =\rho( z) +\({1 \over h}\int_{ x}^{ x+h}{1 \over h}
\int_{ y}^{ y+h}\rho( s-t)\,dt \,ds-\rho( z) \).\nonumber
\end{eqnarray}
By   (\ref{nor14.5}), for $4|h|\leq |z|\leq M$, 
\begin{eqnarray}
&&\Bigg|{1 \over h}\int_{ x}^{ x+h}{1 \over h}
\int_{ y}^{ y+h}\rho( s-t)\,dt \,ds-\rho( z)\Bigg|
\label{ct.2}\\
&&\qquad \le{1 \over h}\int_{0}^{ h}{1 \over h}
\int_{ 0}^{  h}\left|\rho(z+ s-t)-\rho(z)\right|\,dt \,ds   \le  C_M {|h|  \over |z| }\,|\rho( z)|.\nonumber
\end{eqnarray}

Note Ê thatÊ given an integer $j_0$,Ê there exists a $ u_{j_0}>0$, such
that
\be
Ê 1/2\le 1-2j|u|\le (1+u)^j\le 1+2j|u|\le 2
\ee
for all $0\le |u|\le
u_{j_0}$ and $1\le j\le j_0$. Therefore, if we take $ C_{M}{|h|Ê \over |z|
}\le u_{j_0}$ we see that   when $4|h|\le  z\le M$,
\be 
 \({1 \over h}\int_{ x}^{ x+h}{1 \over h}
\int_{ y}^{ y+h}\rho( s-t)\,dt \,ds\)^{j}=   \rho^{j}(z)  +O\(   {|h|Ê \over |z| }\,|\rho^j( z)|\)
\label{ct.3}
\ee
where the   last term  is independent of $z$, (but depends on $M$ and
$j_0$).  

 We estimate the other two integrals in the bracket in  (\ref{nor14.7bk}) similarly to  see that there exists a
constant $C'$ such that for all
$h$ sufficiently smallÊ
 \begin{eqnarray}
Ê Ê\Big| \int\int_{|x-y|\ge C'h} \lc
 \cdots \rc g(x) g( y)\,dx\,dy\Big|&\le& K \int_{C'h\leq |z|\le M}ÊÊ {|h|Ê \over |z|
}\,|\rho^j( z)|\,dz\label{.05}\\
&\le&K'|h|^{1- j\ze}= K' |h|\varphi^j( h).
\nonumber
\end{eqnarray}
where,  in addition to other dependencies, $K$ and $K'$ depend on the
support of
$g$.

By   (\ref{new.2star}) and   the second inequality in (\ref{new.1}), with $z=x-y$,
\be {1 \over h}\int_{ x}^{ x+h}{1 \over h}
\int_{ y}^{ y+h}\rho( s-t)\,dt \,ds\label{est} ={1 \over h^{2}}\int_{ 0}^{  h} 
\int_{ 0}^{  h} \rho(z+u-v  )\,dv \,du .
\ee

We conclude the proof by considering the integral in  
(\ref{nor14.7bk}) when
$|x-y|\le C'h$. Note that 
\bea
{1 \over h} Ê\Big|\int_{ x}^{x+h}\rho( s-y)\,dsÊ\Big|&=&
{1\over h} Ê\Big|\int_{0}^{ h}\rho(z+ s) \,ds Ê\Big|\label{6.23}\\
 &\leÊ &C_{M} {1\over h} \int_{0}^{ h}
{1\over |z+s|^{\ze}}\,ds\nn\\
&\leÊ &2C_{M} {1\over h} \int_{0}^{ h/2}
{1\over | s|^{\ze}}\,ds\nn\\
Ê &\leÊ &
{C\over | h|^{\ze}}=C\varphi(h),\nn
\eea
where for the third  line we use the symmetry of the integrand.
Therefore 
\be
 Ê\Bigg|\({1\over h}\int_{0}^{ h}\rho(z+ s) \,ds\)^j Ê\Bigg|
\le
C\varphi^j(h).
\ee
Similarly 
\be
 Ê\Bigg|\({1 \over h}\int_{ 0}^{ h}{1 \over h}\int_{ 0}^{h}\rho(z+ s-t)\,dt
\,ds\)^{j} Ê\Bigg|\le
C\varphi^j(h).\label{ffwf}
\ee
 Consequently, theÊÊ integral of the first three terms in the bracket in
(\ref{nor14.7bk}), taken over the regionÊ $|z|<C'h$ is bounded by 
$C''h\varphi^j (h).$
As forÊ theÊÊ integral of the last termÊ in the bracket in
(\ref{nor14.7bk}), taken over the regionÊ $|z|<C'h$,
ÊÊ consider 
\be
 Ê\Bigg|\int\!\!\int_ {|x-y|\le C'h} \rho^j(|x-y|) g(x) g( y)\,dx\,dy Ê\Bigg|.
\ee
UsingÊ (\ref{nor14.5a})Ê it is easy to see that this also has the bound 
$C''h\varphi^j (h).$ Thus we obtain (\ref{.03}).  

To obtain  (\ref{.06}) we first note that by (\ref{.m31qz}), a change of variables followed by one integration,     and   (\ref{nor14.5a})
\be
\si^2(h)=2\int_0^h(h-s)\rho(s)\,ds
 \le  Ch^2\varphi(h).\label{.01ja} 
\ee
Therefore
\be
\({\si^2(h)\over h^2}\)^j\le C' \varphi^j(h)\label{.01}.
\ee
The statement in (\ref{.06})
follows immediately from this and (\ref{.03}).\qed

\medskip\noindent{\bf Proof of Corollary \ref{theo-3.2}. } This follows immediately from 
Theorem \ref{theo-3.2n} once we observe that the conditions (\ref{nor14.5a}) and (\ref{nor14.5})
are only used in two places: the proof of Lemma   \ref{calc} which is not need here, and in the proof of (\ref{regularityb}) which is now assumed in condition (\ref{regularity}).\qed

\medskip\noindent{\bf Proof of Corollary \ref{theo-3.1}. } Ê  
 Note
that $ H_{1}(x)=x$. Consequently, for $ f\in L^{ 2}(R^{ 1},\,d\mu)$,  $a_1=
E(\eta f(\eta))$.    Therefore, in Corollary \ref{theo-3.1},
$k_0(f)=1$. 
Hence (\ref{regularity}) is given by the second
condition in (\ref{mm5.17o}).
Also,  By (\ref{si.7a}),
$:(G')^{1}:(I_{[a,b]})=G'(I_{[a,b]})=G(b)-G(a)$.Ê Thus (\ref{mm5.10}) is a special
case ofÊ Ê (\ref{mm5.22}).\qed

\section{FB-mixtures and other examples}\label{sec-5}

Set  \begin{equation}
 \phi (u)=2\int_{-\ff} ^{\ff}  (1-\cos2\pi\la u)\nu(d\la)\label{5.1}
   \end{equation}
   where 
   \begin{equation}
   \int_{-\ff}^{ \ff} (1\wedge\la^{2})\,\nu(d\la)<\ff.\label{5.2}
   \end{equation}
It is well known, (see e.g.  \cite[page 236]{book}),  that $\phi$ can be  the increments variance of Gaussian process with stationary increments that is zero at zero.

   We construct a wide class of examples of Gaussian processes that satisfy the hypotheses of Theorem \ref{theo-3.2n} based on the ideas underlying ``stable mixtures'' considered in \cite[Section 9.6]{book}.
 For $1<\bb<2$ let
  \be
\psi( \la 
)=\int_\bb^2 |\la| ^s\,d\mu(s)\label{q60.2},
\ee  where
$\mu$
  is a 
finite positive measure on $[\bb,2]$ such 
that
\be
\int_\bb^2\frac{d\mu(s)}{2-s}<\ff.\label{08.16}
\ee
  We show in \cite[Section 9.6]{book} that $\psi$ can be represented as in (\ref{5.1}) for some measure $\nu$ satisfying (\ref{5.2}).   Therefore, as we point out in the preceding paragraph,  $\psi$ is the increments variance of a Gaussian process with stationary increments that is zero at zero.
 
  In \cite[Section 9.6]{book} we refer to $\psi$ as a stable mixture because  we were studying  L\'{e}vy processes and $|\la|^{s}$ is the L\'{e}vy exponent   of a symmetric stable process. Here,   since we are concerned with Gaussian processes, we refer to $\psi$ as an FB-mixture because    $|\la|^{s}$ is the increments variance of fractional Brownian motion.

 In \cite[Section 9.6]{book} we study $\psi(\la)$ as $\la\to\ff$. The proofs of 
\cite[Lemma 9.6.1 and Remark 9.6.2]{book}, with obvious modifications, give the proof of the next lemma.

\bl\label{lnormpsi}  The function $\psi(\la)$ is a 
normalized regularly varying
function at zero with index 
$\bb$.   Moreover 
for
$n=1,2,\ldots$,
\be 
\la^n\psi^{(n)}(\la)/\psi(\la) 
\to
\bb(\bb-1)\ldots(\bb-n+1)\qquad\mbox{as}\quad  
\la\to0,\label{09.51a}
\ee
where $\psi^{(n)}$ denotes the $n$--th derivative of $\psi$.
\el

  It follows from (\ref{q60.2}) that $\psi( \la 
)$ is twice differentiable for all $\la\neq 0$ and  
  \be
\psi''( \la 
)=\int_\bb ^2 s(s-1)|\la|^{s-2}\,d\mu(s).\label{q60.2jj}
\ee
We note that  $\psi$ is   convex and  bounded away from the origin. 
In addition, by (\ref{09.51a})
\begin{equation}
   \psi''(\la)\sim \bb(\bb-1) \frac{\psi(\la)}{\la^{2}}\qquad\mbox{as $\la\to 0$}.
   \end{equation}
  It follows   that $\psi''$ is a regularly varying function at zero with index $-(2-\bb)$. Therefore, for any integer $j_{0}\ge 1$ we can find a $1<\bb<2$ such that (\ref{nor14.5a}) holds with $j_{0}\ze<1$. ( Clearly ${1 \over 2}\psi''$ takes the role of $\rho$ in (\ref{rho}).)

 It is easy to see that (\ref{nor14.5}) holds since  
 \begin{equation}
   |\psi''(\la+h)-\psi''(\la)|=\int_\bb ^2 s(s-1)|\la|^{s-2}\(\left|1+\frac{h}{\la}\right|^{s-2}-1\)\,d\mu(s).
   \end{equation}
Lastly, we note that  
\bea
   &&   |\psi(\la+h)+\psi(\la-h)-2\psi (\la)|\\
   &&\qquad=\int_\bb ^2  |\la|^{s } \left|\,\left|1+ h/\la \right|^{s}+\left|1- h/\la \right|^{s}-2\right|  \,d\mu(s)\nn,
   \eea
  which implies (\ref{mm5.9}). Thus we see that   FB-mixtures   $\psi(\la)$ are the increments variance of Gaussian processes that satisfy the hypotheses of Theorem \ref{theo-3.2n}.

    \medskip We give some concrete examples of FB-mixtures.  
    
 \begin{example} \label{ex-5.1}{\rm A simple one that follows immediately from (\ref{q60.2}) is
 \begin{equation}
   \psi(\la)=\sum_{k=0}^{\ff}a_{k}\la^{\bb_{k}}\qquad a_{k}\ge 0,\,\,  \{a_{k}\}\in  \ell_{1},  
   \end{equation}
   where $\bb_{0}=\bb$ and $\{\bb_{k}\}$ is increasing with  $ \sup_{k}\bb_{k}<2$. 
   
 As a slight modification of this, it is easy to see that 
  \begin{equation}
\wt   \psi(\la)=\sum_{k=0}^{n}a_{k}\la^{\bb_{k}}\qquad a_{k}>0,    \end{equation}
where $\bb_{0}=\bb$ and $\{\bb_{k}\}$ is increasing with  $ \bb_{n}=2$, is the increments variance of a Gaussian process, $\wt G$, that satisfies the hypotheses of Theorem \ref{theo-3.2n}.  We get this by taking $\wt G$ to be the sum of two independent Gaussian processes. One with increments variance   the  FB-mixture,   $\psi(\la)=\sum_{k=0}^{n-1}a_{k}\la^{\bb_{k}}$ and the other $\{\sqrt{a_{n}}\, t\, \eta, t\in R_{+}\}$.

 }\end{example}
 
\bl \label{lsmlap}  Let
$\rho(s)$ be a bounded  increasing function on $[0,2-\bb]$,
$1<\bb< 2$ satisfying 
\begin{equation}
   \int_{0}^{2-\bb}\frac{d\rho(v)}{2-\bb-v}<\ff.\label{xx}
   \end{equation}
    Then we can find an FB-mixture with increments variance
\be
\psi(\la)=\la^\bb \hat\rho(\log1/\la)\label{08.14}
   \ee  where
\be
\hat\rho(v)=\int_0^\ff e^{-vs}\,d\rho(s).\label{5.8}
\ee  
\el

\Proof   Let  $\mu(s)$ in (\ref{q60.2}) be a measure 
with 
distribution function
$\rho(s-\bb)$.    It follows from (\ref{xx}) that (\ref{08.16}) holds. Therefore, 
for
$1<\bb<2$,
\bea
\psi(\la)&=&\int_\bb^2\la^s\,d\rho( s-\bb)\\
&=&\la^\bb\int_0^{2-\bb }\la^{s}\,d\rho( 
s)\nn\\
&=&\la^\bb\int_0^{2-\bb }e^{-(\log1/\la)  s}\,d\rho( s)\nn\\ &=& 
\la^\bb
\hat\rho(\log1/\la).\nn
\eea 
\qed

\noindent{\bf  Proof of Proposition \ref{prop-conv}} We use the fact that we can find an FB-mixture of the form (\ref{08.14}). For  $p> 0$, let $\rho(s)\sim s^{p}L(1/s)/\Ga(1+p)$  at zero,  in (\ref{5.8}). Then by \cite[Theorem 14.7.6]{book}, $\hat \rho(s)\sim s^{-p}L(s) $ as $s\to\ff$. Thus (\ref{sam.019}) follows from Lemma \ref{lsmlap}. For $p=0$ if $\rho(s)\sim  L(1/s) $ we must have $   L(1/s) $ increasing as $s$ increases. In this case  $\hat \rho(s)\sim  L(s) $ and $L(s)$ is decreasing. \qed

\begin{example}\label{ex-5.2} {\rm Let \begin{equation}
g(\la)= \frac{\log |\la|-1}{|\la|^{3}}I_{\{|\la|\ge e\}}\label{5.1qq}
   \end{equation}
and consider (\ref{5.1}) with $\nu(d\la)=g(\la)\,d\la$. Then 
\be \si^{2}(u)\approx Cu^{2}\log^{2} 1/u \quad\mbox{  and}\quad  (\si^{2})''(u)\approx \log^{2} 1/u\label{5.18}
\ee
 at zero.  We show this immediately below, and also,  that $\si^{2}(u)$ satisfies the hypotheses of   Theorem \ref{theo-3.2n}. Therefore, we get (\ref{mm5.22nqq}).  
 
\medskip	By (\ref{5.1}) we have 
\bea
 \si^{2} (u)  
 &=&8 \int_{e} ^{\ff}   \sin^{2} \pi\la u  \( \frac{\log \la-1 }{\la^{3} }\) \,d\la \label{5.21}\\
  &\le& 8\pi^{2} u^{2} \int_{e} ^{1/( \pi u)}      \frac{\log \la -1}{\la  }  \,d\la \nn+8\int_{1/( \pi u)} ^{\ff}     \frac{\log \la -1}{\la^{3}  }\,d\la
   \eea
   and
   \begin{equation}
    \si^{2} (u)\ge  Cu^{2} \int_{e} ^{1/( \pi u)}      \frac{\log \la-1 }{\la  }  \,d\la.  \label{5.23a}  \end{equation}
The inequalities in (\ref{5.21}) and (\ref{5.23a}) give the first estimate in (\ref{5.18}).

 Write  $ \sin^{2} \pi\la u=(1-\cos  2\pi\la u)/2$ in (\ref{5.21}) and use trigonometric identities to write  
\bea 
   && { \si^{2} (u+h)-  \si^{2} (u)\over h}\label{5.22}\\
    &&\quad = {4}  \int_{e} ^{\ff} \left\{ {\cos 2\pi\la u (1-\cos 2\pi \la h )+\sin 2\pi\la u\sin 2\pi\la  h\over h}\right\}     \( \frac{\log \la -1}{\la^{3} }\) \,d\la \nn
 \eea 
 Note that the absolute value of the term in the bracket
 \begin{equation}
   \le2\frac{|\sin \pi\la h|+|\sin 2\pi\la h|}{h}\le 6\pi\la.\label{5.23}
   \end{equation}
   Therefore we can use the Dominated Convergence Theorem to see that
 \bea
    (\si^{2} (u))'&=&  8 \pi  \int_{e} ^{\ff}   \sin 2 \pi\la u  \( \frac{\log \la-1 }{\la^{2 } }\) \,d\la.   \label{5.24}
    \eea 
    Using integration by parts we have
    \bea 
      (\si^{2} (u))'  &=&8\pi  \int_{e} ^{\ff}   \( \frac{\log \la -1}{\la^{2 } }\) \,d\int_{0}^{\la} \sin 2 \pi s  u \,ds\label{wew}\\
    &=&{4\over u}  \int_{e} ^{\ff}  (1-\cos 2 \pi\la u)  \( \frac{2\log \la -3}{\la^{3 } }\) \,d\la  \nn .
    \eea
  Exactly the same argument  used in (\ref{5.22}) and (\ref{5.23})   shows that    
    \bea
    (\si^{2} (u))''&=&-{4\over u^{2}}  \int_{e} ^{\ff}  (1-\cos2 \pi\la u)  \( \frac{2\log \la -3}{\la^{3 } }\)\,d\la  \label{5.25}  \\
    & &\qquad+  {8\pi\over u }   \int_{e} ^{\ff}   \sin 2\pi\la u \( \frac{2\log \la-3 }{\la^{2 } }\) \,d\la \nn\\
    &=&I+II\nn.
    \eea
    Considering (\ref{5.21}) we see that 
    \begin{equation}
   |I|\approx C'\frac{  \si^{2} (u)}{u ^{2}}\label{5.27},
   \end{equation}
   and by the same methods used to obtain the first estimate in (\ref{5.18}) we get
   \bea 
   II\approx (\log 1/u)^{2}.
 \eea 
Thus we get the second estimate in (\ref{5.18}).
   
   It follows from (\ref{5.21}) that 
   \bea
    \si^{2} (au)&=&8a^{2} \int_{ea} ^{\ff}   \sin^{2} \pi\la u  \( \frac{\log (\la/a) -1}{\la^{3} }\) \,d\la \label{5.28}  \\ 
   &= & a^{2}\( \si^{2}(u)-8\int_{e } ^{ea}   \sin^{2} \pi\la u  \( \frac{\log  \la  -1}{\la^{3} }\) \,d\la \right.  \nn\\
   &&\left.\qquad +8\log 1/a\int_{e a} ^{\ff} {  \sin^{2} \pi\la u \over \ \la^{3} }  \,d\la \)\nn.
    \eea 
It is easy   to see that the last two integrals in (\ref{5.28})  are $o( \si^{2} (u))$ at zero. This shows that $\si^{2}(u)$ is regularly varying at zero with index 2. 

 Note that by (\ref{5.18})
  \begin{equation}  
    (\si^{2})'' (s)=O\(\frac{  \si^{2}(s)}{s^{2}}\)
  \end{equation}
  at zero. This implies that $\si^{2}$ satisfies (\ref{mm5.9}). 
  
  We now show that $\si^{2}$ satisfies  (\ref{nor14.5}).  We write the integral II in (\ref{5.25}) as
  \begin{equation}
    {8\pi\over u }  \( \int_{e} ^{\ff}   \sin 2\pi\la u \( \frac{2\log \la-2 }{\la^{2 } }\) \,d\la-\int_{e} ^{\ff}  { \sin 2\pi\la u\over   \la^{2 } }  \,d\la\).
   \end{equation}
  Clearly
  \begin{equation}
Q(u)=  {8\pi\over u } \int_{e} ^{\ff}  { \sin 2\pi\la u\over   \la^{2 } }  \,d\la=16\pi^{2}\int_{2\pi u e}^{\ff}\frac{\sin s}{s^{2}}\,ds.
   \end{equation}
  Furthermore, by integration by parts, as in (\ref{wew}),
  \bea 
 &&  {8\pi\over u }   \int_{e} ^{\ff}   \sin 2\pi\la u \( \frac{2\log \la-2 }{\la^{2 } }\) \,d\la\\
 &&\qquad={8\over u^{2}}  \int_{e} ^{\ff}  (1-\cos2 \pi\la u)  \( \frac{2\log \la -3}{\la^{3 } }\)\,d\la\nn.
 \eea
  Substituting this in (\ref{5.25}) we see that 
   \be 
    (\si^{2} (u))''= {4\over u^{2}}  \int_{e} ^{\ff}  (1-\cos2 \pi\la u)  \( \frac{2\log \la -3}{\la^{3 } }\)\,d\la+Q(u).  
    \ee  
 As in (\ref{wew}) we can differentiate under the integral sign to get
 \bea 
      (\si^{2} (u))'''&= &-{8\over u^{3}}  \int_{e} ^{\ff}  (1-\cos2 \pi\la u)  \( \frac{2\log \la -3}{\la^{3 } }\)\,d\la\\
      &&\quad+{8\pi\over u^{2}}  \int_{e} ^{\ff}  \sin 2 \pi\la u \( \frac{2\log \la -3}{\la^{2 } }\)\,d\la+ O(1/u).\nn 
 \eea
  Separating the integral as in (\ref{5.21}) we see that 
  \begin{equation}
   |   (\si^{2} (u))'''|\le C\frac{(\log 1/u)^{2}}{u}.
   \end{equation}
  This implies    that $\si^{2}$ satisfies  (\ref{nor14.5}).
  Thus the Gaussian process determined   by (\ref{5.1}) and (\ref{5.1qq}) satisfies the  hypotheses of   Theorem \ref{theo-3.2n}.}\end{example}

\def\noopsort#1{} \def\printfirst#1#2{#1}
\def\singleletter#1{#1}
             \def\switchargs#1#2{#2#1}
\def\bibsameauth{\leavevmode\vrule height .1ex
             depth 0pt width 2.3em\relax\,}
\makeatletter
\renewcommand{\@biblabel}[1]{\hfill#1.}\makeatother
\newcommand{\bysame}{\leavevmode\hbox to3em{\hrulefill}\,}

\end{document}